\documentclass[leqno]{revtex4}
\usepackage[utf8]{inputenc}
\usepackage[russian]{babel}
\usepackage{amssymb,amsmath,amsthm}
\renewcommand{\subsection}{\refstepcounter{subsection}%
\par\bigskip\noindent\textbf{\upshape\arabic{subsection}. }}
\renewcommand{\subsubsection}{\refstepcounter{subsubsection}%
\par\medskip\noindent\textbf{\upshape\arabic{subsection}.%
\arabic{subsubsection}.\ }}

\numberwithin{equation}{section}

\DeclareMathOperator{\ind}{ind}
\newcommand{\Wo}{{\raisebox{0.2ex}{\(\stackrel{\circ}{W}\)}}{}}

\begin{document}
\renewcommand{\proofname}{\upshape Д\:о\:к\:а\:з\:а\:т\:е\:л\:ь\:с\:т\:в\:о}
\title{О вычислении собственных значений задачи Штурма--Лиувилля
с индефинитным весом, имеющим самоподобную первообразную}
\author{А.~А.~Владимиров}
\email[Электронный адрес: ]{vladimi@mech.math.msu.su}
\affiliation{Московский Государственный Университет~им.~М.~В.~Ломоносова,\\
механико-математический факультет}
\thanks{Работа поддержана РФФИ, грант \mbox{\No~04-01-00712}, грантом
поддержки ведущих научных школ \mbox{НШ-5247.2006.1}, и INTAS, грант
\mbox{\No\,05-1000008-7883}.}

\noindent УДК~517.984.5+519.614.2+519.671
\begin{abstract}
В статье описывается эффективный метод вычисления собственных значений
граничной задачи
\begin{gather*}
	-y''-\lambda\rho y=0,\\ y(0)=y(1)=0,
\end{gather*}
где \(\rho\) есть функция из пространства \(\Wo_2^{-1}[0,1]\), имеющая
самоподобную обобщённую первообразную \(P\in L_2[0,1]\).
\end{abstract}

\maketitle

\section*{Введение}
\subsection
Целью настоящей статьи является изложение способа вычисления точных оценок
собственных значений граничной задачи
\begin{gather}\label{eq:0:1}
	-y''-\lambda\rho y=0,\\ \label{eq:0:2} y(0)=y(1)=0,
\end{gather}
где весовая функция \(\rho\) принадлежит пространству \(\Wo_2^{-1}[0,1]\)
и имеет самоподобную обобщённую первообразную \(P\in L_2[0,1]\). Этот метод
был использован при получении численных примеров, приведённых
в работах~\cite{VlaSh1} и~\cite{VlaSh2}.

\section{Сведение задачи к конечномерной}
\subsection
Как и в работах~\cite{VlaSh1} и~\cite{VlaSh2}, под \emph{индексом инерции}
\(\ind F\) эрмитова оператора \(F\), действующего в некотором гильбертовом
пространстве \(\mathfrak H\), мы будем понимать точную верхнюю грань
размерностей конечномерных подпространств \(\mathfrak M\subseteq\mathfrak H\),
удовлетворяющих условию
\begin{equation}\label{eq:2:2}
	\exists\varepsilon>0\,\forall x\in\mathfrak M\qquad
	\langle Fx,x\rangle\leqslant -\varepsilon\,\|x\|^2_{\mathfrak H}.
\end{equation}
Имеет место следующий факт:

\subsubsection
{\itshape Пусть \(\mathfrak H_1\) и \(\mathfrak H_2\) "--- два гильбертовых
пространства, а \(\mathfrak H\) "--- прямая сумма \(\mathfrak H_1\oplus
\mathfrak H_2\) пространств \(\mathfrak H_1\) и \(\mathfrak H_2\). Пусть также
\(F:\mathfrak H\to\mathfrak H\) "--- ограниченный эрмитов оператор
с блочно-матричным представлением
\[
	F=\begin{pmatrix}A&B^*\\B&C\end{pmatrix},
\]
в котором оператор \(C:\mathfrak H_2\to\mathfrak H_2\) равномерно положителен.
Тогда выполняется равенство
\begin{equation}\label{eq:2:1}
	\ind F=\ind (A-B^*C^{-1}B).
\end{equation}
}

\begin{proof}
Ввиду ограниченной обратимости оператора \(C\), для оператора \(F\) можно
рассмотреть факторизацию Фробениуса--Шура \(F=H^*F_0H\), где через \(H\)
и \(F_0\) обозначены операторы с блочно-матричными представлениями
\begin{align*}
	H&\rightleftharpoons\begin{pmatrix}1&0\\ C^{-1}B&1\end{pmatrix},&
	F_0&\rightleftharpoons\begin{pmatrix}A-B^*C^{-1}B&0\\ 0&C\end{pmatrix}.
\end{align*}
При этом оператор \(H\) обладает ограниченным обратным вида
\begin{align*}
	H^{-1}&=\begin{pmatrix}1&0\\ -C^{-1}B&1\end{pmatrix},
\end{align*}
а для любого вектора \(x\in\mathfrak H\) выполняются равенства
\begin{flalign*}
	&& \langle Fx,x\rangle&=\langle H^*F_0H\,x,x\rangle\\
	&& &=\langle F_0\,(Hx),(Hx)\rangle.&&
\end{flalign*}
Тогда для любого конечномерного подпространства \(\mathfrak M\subseteq
\mathfrak H\), удовлетворяющего условию~\eqref{eq:2:2}, подпространство
\(\mathfrak M_0\rightleftharpoons H\,\mathfrak M\) имеет ту же размерность
и удовлетворяет условию
\begin{equation}\label{eq:2:3}
	\exists\varepsilon>0\,\forall x\in\mathfrak M_0\qquad
	\langle F_0x,x\rangle\leqslant -\varepsilon\,\|x\|^2_{\mathfrak H},
\end{equation}
а для любого конечномерного подпространства \(\mathfrak M_0\subseteq
\mathfrak H\), удовлетворяющего условию~\eqref{eq:2:3}, подпространство
\(\mathfrak M\rightleftharpoons H^{-1}\,\mathfrak M_0\) имеет ту же
размерность и удовлетворяет условию~\eqref{eq:2:2}. Таким образом, выполняется
равенство \(\ind F=\ind F_0\). Ввиду положительности оператора \(C\),
это означает выполнение равенства~\eqref{eq:2:1}.
\end{proof}

\subsection
В дальнейшем через \(\mathfrak H\) мы будем обозначать пространство Соболева
\(\Wo_2^1[0,1]\), снабжённое скалярным произведением
\[
	\langle y,z\rangle\rightleftharpoons\int\limits_0^1 y'
	\overline{z'}\,d\mu,
\]
где \(d\mu\) "--- линейная мера Лебега. Простым следствием
теоремы~\cite[Теорема~4.1]{VlaSh1} является такой факт:

\subsubsection
{\itshape Пусть \(F\) "--- пучок действующих в пространстве \(\mathfrak H\)
линейных операторов, удовлетворяющий тождеству
\begin{equation}\label{eq:2:4}
	\forall\lambda\in\mathbb R\;\forall y\in\mathfrak H\qquad
	\langle F(\lambda)y,y\rangle=\int\limits_0^1\biggl(|y'|^2+\lambda P\cdot
	(|y|^2)'\biggr)\,d\mu.
\end{equation}
Пусть также \(\{\nu_n\}_{n=1}^{\infty}\) "--- последовательность (возможно,
частичная) сосчитанных в порядке возрастания положительных собственных
значений граничной задачи~\eqref{eq:0:1},~\eqref{eq:0:2}. Тогда для любых
натурального числа \(n\geqslant 1\) и вещественного числа \(\lambda\in
(0,\nu_n)\) выполняется неравенство \(\ind F(\lambda)<n\), а для любых
натурального числа \(n\geqslant 1\) и вещественного числа \(\lambda\in
(\nu_n,+\infty)\) выполняется неравенство \(\ind F(\lambda)\geqslant n\).
}

\bigskip
Таким образом, при наличии в нашем распоряжении метода вычисления достаточно
точных оценок индексов инерции операторов из пучка \(F\) мы можем вычислять
оценки собственных значений граничной задачи~\eqref{eq:0:1},~\eqref{eq:0:2}
на основе метода деления отрезка.

\subsection
Пусть теперь \(S\) "--- набор из натурального числа \(N>1\) и вещественных
чисел \(a_k>0\), \(d_k\) и \(\beta_k\), где \(k=1,\ldots,N\), удовлетворяющих
соотношениям
\begin{align*}
	\sum_{k=1}^N a_k&=1,&\theta_S\rightleftharpoons
	\sqrt{\sum\limits_{k=1}^N a\,|d_k|^2}&<1.
\end{align*}
Как и в работе~\cite{VlaSh1}, с набором \(S\) мы связываем семейство
\(\{G_k\}_{k=1}^{N}\) действующих в пространстве \(L_2[0,1]\) линейных
операторов вида
\[
	\forall (\zeta,\xi)\subseteq [0,1]\qquad
	G_k\,\chi_{(\zeta,\xi)}=\chi_{(\alpha_k+a_k\zeta,\alpha_k+a_k\xi)},
\]
где через \(\{\alpha_k\}_{k=1}^{N+1}\) обозначен набор чисел вида
\(\alpha_1=0\), \(\alpha_{k+1}=\alpha_k+a_k\), а через \(\chi_I\) "---
индикаторы промежутков \(I\). При этом на основе семейства операторов
\(\{G_k\}_{k=1}^N\) мы конструируем нелинейный, вообще говоря, оператор
\(G_S\) вида
\begin{equation}\label{eq:2:6}
	\forall f\in L_2[0,1]\qquad G_S\,f\rightleftharpoons
	\sum\limits_{k=1}^N (d_k\cdot G_kf+\beta_k\cdot\chi_{(\alpha_k,
	\alpha_{k+1})}).
\end{equation}
Этот оператор является сжимающим, и введённая ранее величина \(\theta_S\)
представляет собой его коэффициент сжатия (см.~\cite[Лемма~3.1]{VlaSh1}).

Сжимающие операторы, допускающие представление в виде~\eqref{eq:2:6},
мы называем \emph{операторами подобия}. Функции \(f\in L_2[0,1]\), являющиеся
неподвижными точками операторов подобия, мы называем \emph{самоподобными}.
Набор чисел \(S\), которому отвечает оператор подобия \(G_S\), оставляющий
неподвижной некоторую заранее фиксированную функцию \(f\in L_2[0,1]\),
мы называем набором \emph{параметров самоподобия} функции \(f\).

\subsection
С произвольно фиксированным набором \(S\) параметров самоподобия функции
\(P\) мы в дальнейшем будем связывать конечномерное подпространство
\(\mathfrak H_{S,1}\) пространства \(\mathfrak H\), обладающее базисом
\(\{y_{S,k}\}_{k=1}^{N-1}\) вида
\[
	y_{S,k}(x)=\left\{\begin{array}{ll}\dfrac{x-\alpha_k}{a_k}&
	\text{при } x\in [\alpha_k,\alpha_{k+1}],\\
	\dfrac{\alpha_{k+2}-x}{a_{k+1}}&
	\text{при } x\in [\alpha_{k+1},\alpha_{k+2}],\\
	0&\text{иначе.}\end{array}\right.
\]
При этом через \(\mathfrak H_{S,2}\) мы будем обозначать ортогональное
дополнение \(\mathfrak H\ominus\mathfrak H_{S,1}\) подпространства
\(\mathfrak H_{S,1}\). Нетрудно видеть, что справедливо следующее утверждение:

\subsubsection
{\itshape Подпространство \(\mathfrak H_{S,2}\) имеет вид \(\mathfrak H_{S,2}=
\{y\in\mathfrak H\mid\forall k=2,\ldots,N\quad y(\alpha_k)=0\}\).
}

\bigskip
Введём теперь в рассмотрение три пучка \(A_S\), \(B_S\) и \(C_S\) линейных
операторов, значениями которых \(A_S(\lambda):\mathfrak H_{S,1}\to
\mathfrak H_{S,1}\), \(B_S(\lambda):\mathfrak H_{S,1}\to\mathfrak H_{S,2}\)
и \(C_S(\lambda):\mathfrak H_{S,2}\to\mathfrak H_{S,2}\) при произвольно
фиксированном \(\lambda\in\mathbb R\) являются элементы блочно-матричного
представления
\[
	F(\lambda)=\begin{pmatrix}A_S(\lambda)&B_S^*(\lambda)\\
	B_S(\lambda)&C_S(\lambda)\end{pmatrix}
\]
определённого условием~\eqref{eq:2:4} оператора \(F(\lambda)\). Имеет место
следующий факт:

\subsubsection
{\itshape Пусть даны два вещественных числа \(\lambda>0\) и \(\varepsilon>0\),
удовлетворяющие неравенствам \(\varepsilon\geqslant 2\|B_S(\lambda)\|^2\)
и \(\|C_S(\lambda)-1\|<1/2\). Тогда выполняются неравенства
\(\ind A_S(\lambda)\leqslant\ind F(\lambda)\) и \(\ind F(\lambda)\leqslant
\ind (A_S(\lambda)-\varepsilon)\)} [1.1].

\bigskip
Утверждение~4.2 позволяет свести задачу вычисления оценок индекса
инерции оператора \(F(\lambda)\) к допускающей непосредственное решение
на ЭВМ задаче вычисления оценок индексов инерции конечномерных операторов
\(A_S(\lambda)\) и \(A_S(\lambda)-\varepsilon\). Однако при этом встаёт вопрос
об области применимости утверждения~4.2 и о степени точности
получаемых на его основе оценок величины \(\ind F(\lambda)\). Изучению
этого вопроса будет посвящён следующий параграф.

\section{Исследование точности конечномерных приближений}
\subsection
Имеет место следующий факт:

\subsubsection
{\itshape Пусть даны набор \(S\) параметров самоподобия функции \(P\)
и вещественное число \(\lambda>0\). Тогда выполняются неравенства
\(\|B_S(\lambda)\|\leqslant\lambda\theta_S\,\|P\|_{L_2[0,1]}\)
и \(\|C_S(\lambda)-1\|\leqslant\lambda\theta_S\,\|P\|_{L_2[0,1]}\).
}

\begin{proof}
Заметим, что полуторалинейные формы операторов \(B_S(\lambda)\)
и \(C_S(\lambda)-1\) удовлетворяют тождествам
\begin{gather*}
	\forall y\in\mathfrak H_{S,1}\;\forall z\in\mathfrak H_{S,2}\qquad
	\langle B_S(\lambda)y,z\rangle=\lambda\cdot\int\limits_0^1
	P\cdot(y\overline{z})'\,d\mu,\\
	\forall y,z\in\mathfrak H_{S,2}\qquad
	\langle (C_S(\lambda)-1)y,z\rangle=\lambda\cdot\int\limits_0^1
	P\cdot(y\overline{z})'\,d\mu.
\end{gather*}
При этом функция \(P\) в правых частях выписанных тождеств может быть заменена
функцией \(\hat P_S\in L_2[0,1]\) вида \(\hat P_S\rightleftharpoons
\sum_{k=1}^N d_k\cdot G_kP\) [\S~1.4.1], удовлетворяющей
очевидному равенству \(\|\hat P_S\|_{L_2[0,1]}=\theta_S\,\|P\|_{L_2[0,1]}\).
Заметим также, что для любых двух функций \(y,z\in\mathfrak H\) выполняются
соотношения
\begin{flalign*}
	&& \|(\overline{y}z)'\|_{L_2[0,1]}&\leqslant
	\|\overline{y}z'\|_{L_2[0,1]}+\|\overline{y}'z\|_{L_2[0,1]}\\
	&& &\leqslant \|y\|_{C[0,1]}\cdot \|z'\|_{L_2[0,1]}+
	\|z\|_{C[0,1]}\cdot \|y'\|_{L_2[0,1]}\\
	&& &\leqslant \dfrac{\|y'\|_{L_2[0,1]}}{2}\cdot \|z'\|_{L_2[0,1]}+
	\dfrac{\|z'\|_{L_2[0,1]}}{2}\cdot \|y'\|_{L_2[0,1]}\\
	&& &=\|y\|_{\mathfrak H}\cdot \|z\|_{\mathfrak H}.&&
\end{flalign*}
С учётом сделанных замечаний, доказываемое утверждение тривиальным образом
выводится из неравенства Коши--Буняковского.
\end{proof}

\subsection
Имеют место следующие три простых факта:

\subsubsection
{\itshape Пусть даны набор \(S\) параметров самоподобия функции \(P\)
и натуральное число \(m>1\). Тогда оператор \(G_S^m\) является оператором
подобия с неподвижной точкой \(P\).
}

\subsubsection
{\itshape Пусть дано вещественное число \(\varepsilon>0\). Тогда существует
набор \(S\) параметров самоподобия функции \(P\), удовлетворяющий неравенству
\(\theta_S<\varepsilon\)} [2.1].

\subsubsection
{\itshape Пусть даны набор \(S\) параметров самоподобия функции \(P\)
и два вещественных числа \(\lambda>0\) и \(\varepsilon>0\), удовлетворяющие
неравенствам \(\varepsilon\geqslant 2\lambda^2\theta_S^2\,\|P\|^2_{L_2[0,1]}\)
и \(\lambda\theta_S\,\|P\|_{L_2[0,1]}<1/2\). Тогда выполняются неравенства
\(\ind A_S(\lambda)\leqslant\ind F(\lambda)\) и \(\ind F(\lambda)\leqslant
\ind(A_S(\lambda)-\varepsilon)\)} [\S~1.4.2, 1.1].

\bigskip
Утверждения~2.2 и~2.3 указывают способ вычисления
оценок индекса инерции оператора \(F(\lambda)\) при произвольно фиксированном
значении \(\lambda>0\).

\subsection
Пусть теперь \(\nu_n\) "--- имеющее произвольно фиксированный номер
\(n\geqslant 1\) положительное собственное значение
задачи~\eqref{eq:0:1},~\eqref{eq:0:2}.

Заметим, что при любом значении \(\lambda\in (0,\nu_n)\) существует
вещественное число \(\delta\in (0,1)\), удовлетворяющее неравенству \(\lambda<
(1-\delta)\nu_n\). В таком случае для любых набора \(S\) параметров
самоподобия функции \(P\) и вещественного числа \(\varepsilon\in (0,\delta]\),
удовлетворяющих условиям утверждения~2.3, выполняются соотношения
\begin{flalign*}
	&& \ind (A_S(\lambda)-\varepsilon)&\leqslant
	\ind (A_S(\lambda)-\delta)\\
	&& &\leqslant\ind (F(\lambda)-\delta)\\
	&& &=\ind F\left(\dfrac{\lambda}{1-\delta}\right)&
	\text{[\eqref{eq:2:4}]}&\\
	&& &<n.&\text{[\S~1.2.1]}&
\end{flalign*}

Заметим также, что при любом значении \(\lambda\in (\nu_n,+\infty)\) оператор
\(F(\lambda)-1\) является вполне непрерывным [2.2, 1.1]. Поэтому существует
вещественное число \(\delta\in (0,1/2)\), удовлетворяющее равенству
\(\ind (F(\lambda)+\delta)=\ind F(\lambda)\) (см.~\cite[пункт~95]{RN}).
В таком случае для любого набора \(S\) параметров самоподобия функции \(P\),
удовлетворяющего неравенству \(2\lambda^2\theta_S^2\,
\|P\|^2_{L_2[0,1]}\leqslant\delta\), выполняются соотношения
\begin{flalign*}
	&& \ind A_S(\lambda)&\geqslant\ind (A_S(\lambda)+\delta-B_S(\lambda)^*
	C_S^{-1}(\lambda)B_S(\lambda))&
	\text{[1.1]}&\\
	&& &\geqslant\ind (F(\lambda)+\delta)&
	\text{[\S~1.1.1]}&\\
	&& &=\ind F(\lambda)\\
	&& &\geqslant n.&\text{[\S~1.2.1]}&
\end{flalign*}

Таким образом, получаемые на основе утверждений~2.2 и~2.3 оценки величин
\(\ind F(\lambda)\) позволяют для любых двух вещественных чисел
\(\lambda_1>0\) и \(\lambda_2>\lambda_1\) установить
верность одного из неравенств \(\lambda_1\leqslant\nu_n\)
или \(\lambda_2\geqslant\nu_n\). Иначе говоря, они позволяют находить
при помощи метода деления отрезка сколь угодно точные оценки собственного
значения \(\nu_n\).

\section{Непосредственное вычисление оценок}
\subsection
Использование условия самоподобия функции \(P\) позволяет выписывать
рекуррентные формулы для её степенных моментов
\[
	\mathbf P_q\rightleftharpoons\int\limits_0^1 P\cdot x^q\,d\mu,\qquad
	q=0,1,2,\ldots
\]
В частности, употребляемые в дальнейшем моменты нулевой и первой степеней
имеют вид
\begin{align*}
	\mathbf P_0&=\dfrac{\sum\limits_{k=1}^N a_k\,\beta_k}{1-
	\sum\limits_{k=1}^N a_k\,d_k},&
	\mathbf P_1&=\dfrac{\sum\limits_{k=1}^N a_k\cdot\left(
	\dfrac{a_k\,\beta_k}{2}+\alpha_k\,d_k\cdot \mathbf P_0+
	\alpha_k\,\beta_k\right)}{1-\sum\limits_{k=1}^N a_k^2\,d_k}.
\end{align*}

\subsection
Прямым просчётом устанавливается следующий факт:

\subsubsection
{\itshape Пусть даны набор \(S\) параметров самоподобия функции \(P\)
и два вещественных числа \(\lambda>0\) и \(\varepsilon\geqslant 0\). Тогда
матрица квадратичной формы оператора \(A_S(\lambda)-\varepsilon\) в базисе
\(\{y_{S,k}\}_{k=1}^{N-1}\) является трёхдиагональной и имеет элементы
\begin{align*}
	\langle (A_S(\lambda)-\varepsilon)\, y_{S,k},y_{S,k}\rangle&=
	(1-\varepsilon)\cdot [a_k^{-1}+a_{k+1}^{-1}]+
	\lambda\cdot\left[2d_{k+1}\,(\mathbf P_1-\mathbf P_0)+2d_k\,\mathbf P_1
	+\beta_k-\beta_{k+1}\right],\\
	\langle (A_S(\lambda)-\varepsilon)\, y_{S,k},y_{S,k-1}\rangle&=
	-(1-\varepsilon)\cdot a_k^{-1}-\lambda\cdot d_k\,
	(2\mathbf P_1-\mathbf P_0).
\end{align*}
}

Сигнатура матрицы из утверждения~2.1 может теперь быть вычислена
различными способами "--- например, как число перемен знака в ряде главных
миноров этой матрицы.

\subsection
Машинная программа на языке РЕФАЛ, вычисляющая оценки положительных
собственных значений задачи~\eqref{eq:0:1},~\eqref{eq:0:2} на основе
изложенной схемы, может быть найдена по адресу
\texttt{http://www.math.msu.su/labs/spectrallab/soft/Dirichlet.ref}.

\end{document}